\documentclass[12pt]{amsart}
\usepackage{graphicx}
\usepackage{amssymb}
\usepackage{float}
\usepackage{placeins}
\usepackage[utf8]{inputenc}

\newtheorem{theorem}{Theorem}[section]

\newtheorem{corollary}[theorem]{Corollary}
\newtheorem{definition}{Definition}

\author{A. Kanel-Belov}
\address{\textbf{Alexei Kanel-Belov:} Department of Mathematics, Bar-Ilan university, Ramat Gan 52100, Israel.}
\email{kanelster@gmail.com}
\thanks{The first author supported by grant RNF 22-11-00177}

\author{M. Golafshan}
\address{\textbf{Mehdi Golafshan:} Department of Discrete Mathematics, Moscow institute of physics and technology,  Dolgoprudny 141700, Russia.}
\email{m.golafshan@phystech.edu}
\thanks{}

\author{S. Malev}
\address{\textbf{Sergey Malev:} Department of Mathematics, Ariel university, Ariel 40700, Israel.}
\email{sergeyma@ariel.ac.il}
\thanks{}

\author{R. Yavich}
\address{\textbf{Roman Yavich:}  Department of Mathematics, Ariel university, Ariel 40700, Israel.}
\email{romany@ariel.ac.il}
\thanks{}

\thanks{We would like to thank A.Domoshnitsky for interesting and fruitful
discussions regarding this paper.}
\title[area and perimeter full distribution
functions]{Area and perimeter full distribution  functions for planar Poisson line processes and Voronoi diagrams}

\date{}
\begin{document}

\begin{abstract}
The challenges of examining random partitions of space are a significant class of problems in the theory of geometric transformations.
Richard Miles calculated moments of areas and perimeters of any order (including expectation) of the random division of space in 1972.
In the paper we calculate whole distribution function of random divisions of plane by Poisson line process. Our idea is to interpret a random polygon as the evolution of a segment along a moving straight line. In the plane example, the issue connected with an infinite number of parameters is overcome by considering a secant line.

We shall take into account the following tasks:

\textbf{1.}
On the plane, a random set of straight lines is provided, all shifts are equally likely, and the distribution law is of the form $F(\varphi).$
What is the area distribution of the partition's components?

 \textbf{2.}
On the plane, a random set of points is marked.
Each point $A$ has an associated area of attraction, which is the collection of points in the plane to which the point $A$ is the nearest of the designated ones.

In the first problem, the density of moved sections adjacent to the line allows for the expression of the balancing ratio in kinetic form.
Similarly, one can write the perimeters’ kinetic equations. We will demonstrate how to reduce these equations to the Riccati equation using the Laplace transformation in this paper.

In fact, we formulate the distribution function of area and perimeter and the joint distribution of them with a Poisson line process based on differential equations. Also, for Voronoi diagrams. These are the main search results (see theorems \ref{thm1}, \ref{thm2}, \ref{thm3}).

\end{abstract}

\maketitle

\noindent\rule{12.7cm}{1.0pt}



\noindent\rule{12.7cm}{1.0pt}

\section{Introduction}\label{intro}

Statistical geometry or stochastic geometry is an interdisciplinary field between pure and applied mathematics.
On the other hand, due to its many applications with other branches, it can be called multidisciplinary.
Among its most important uses are astronomy, telecommunications,  image analysis and stereology.
In his PhD dissertation, Alexei Kanel Belov worked on statistical geometry and its applications \cite{13,11}.

An important class of problems in the theory of geometric transformations are the problems of studying random partitions of space.
Similar problems arise from mining. In the general case, as a rule, it is possible to find only moments of random variables associated with the partition, the problems arising here are very sophisticated.

Numerous models for the homogeneous random partition of three-dimensional space are surveyed.
Richard Miles discussed about the partitioning of space by random planes, followed by the Voronoi tessellation and the dual Delaunay tessellation of tetrahedra, both of which are formed from random points in space \cite{7}.
 Also in 1969, he talked about the distributions of the typical polygon and of the fundamental polygon are
fairly different: $(\mathrm{i})$ hose of the typical polygon are weighted in number;
$(\mathrm{ii})$ hose of the fundamental polygon are weighted in area \cite{miles}.

According to the \textit{zero-one law}, there are questions about how the justification works. The partition is made up of an infinite set of random variables, and all properties that do not depend on the value of a small number of these variables are met with a probability of $0$ or $1$. Such features are asymptotic in nature. To substantiate, the concept of approximation to the case of a continuous distribution of $n$ system cases is quite beneficial.

In the plane situation, a partial differential equation can be used to explain the distribution of a number of quantities associated with a random partition, which can be simplified to the \textit{Riccati equation} in some cases.

The reason for this is that in the plane example, the problem connected with the fact that the partition's components have an infinite number of parameters is overcome by considering a secant line \cite{14,15}. The idea is that random partitioning of a plane can be thought of as a process along a moving straight line with few parameters (section line, arrea and perimeter passed etc.). The advantage of this strategy is that the polygon's section is a segment, and for future \textit{evolution}, just its length and the velocities of its left and right \textit{ends} are significant. (When studying the law of distribution of systems, areas, or perimeters, we must \textit{remember} the previously moved region (perimeter)).

For space, the concept of a moving secant is works not so nice: a polyhedron's section might be any polygon with an arbitrary number of sides specified by the corresponding number of parameters. However some equations can be written. Note that density of polyhedra type exponentially decrease respect to number of sides.

Local independence of moving its ends can be used for simplifying equations (see theorem \ref{thm2} and section \ref{ScEndsIndep}).
Finally we came from distribution among section to usual distribution using parameter $t$ (projection on the line orthogonal to the section line $L$ (see theorem \ref{thm3} and section \ref{SctoUsualDistr}).

In forthcoming paper we suppose to use similar ideas for nonelection plane and sphere.


The last section will discuss the justification.

\subsection{Main results}
These result is obtained by line Poisson process.
Here, $t$ is the projection of polygon to the perpendicular of line.
In this paper we will prove following theorem:

\begin{theorem}{\label{thm1}}
The distribution function along the lines can be represented by kinetic equations as follow:
\begin{multline*}
l\frac{\partial N(S,l,\alpha_1,\alpha_2,t)}{\partial S} + \left(\cot (\alpha_1) +\cot (\alpha_2)\right)\frac{\partial N(S,l,\alpha_1,\alpha_2,t)}{\partial l} +\frac{\partial N(S,l,\alpha_1,\alpha_2,t)}{\partial t} + \\ \lambda N(S,l,\alpha_1,\alpha_2,t) \left( \tan (\frac{\alpha_1}{2}) + \tan (\frac{\alpha_2}{2})\right )-
\lambda  \cdot \int\limits_0^{\alpha_1} \frac{N(S,l,\varphi,\alpha_2,t) \sin (\alpha_1-\varphi)}{\sin(\varphi)} \mathrm{d} \varphi - \\
\lambda  \cdot \int\limits_0^{\alpha_2} \frac{N(S,l,\alpha_1,\varphi,t) \sin (\alpha_2-\varphi)}{\sin(\varphi)} \mathrm{d} \varphi =0.
\end{multline*}

Similarly,  the kinetic equation of the perimeters:
\begin{multline*}
\left(\frac{1}{\cos (\alpha_1)} + \frac{1}{\cos (\alpha_2)}\right)\frac{\partial N}{\partial P} + \frac{\partial N(P,l,\alpha_1,\alpha_2,t)}{\partial t}
+\frac{\partial N(P,l,\alpha_1,\alpha_2,t)}{\partial l}\left(\cot (\alpha_1) +\cot (\alpha_2)\right)+\\ \lambda N(P,l,\alpha_1,\alpha_2,t) \left( \tan(\frac{\alpha_1}{2}) + \tan(\frac{\alpha_2}{2})  \right) -
\lambda  \cdot \int\limits_0^{\alpha_1} \frac{N(P,l,\varphi,\alpha_2,t) \sin (\alpha_1-\varphi)}{\sin(\varphi)} \mathrm{d} \varphi -\\
\lambda  \cdot \int\limits_0^{\alpha_2} \frac{N(P,l,\alpha_1,\varphi,t) \sin (\alpha_2-\varphi)}{\sin(\varphi)} \mathrm{d} \varphi =0.
\end{multline*}

Let $\mathrm{d}F$ be the expectation of the number of intersection points of a unit segment normal to a straight line with an angle of deviation $\varphi$ with a straight line whose angle of deviation is enclosed between $\varphi$ and $\varphi+ \mathrm{d} \varphi$.
We define the function $F(\varphi)$ as the limit $F(\varphi)=\lim_{\mathrm{d} \varphi\rightarrow 0}\frac{\mathrm{d} F}{\mathrm{d}\varphi}$.
The kinetic equation of areas can be written as:

\begin{multline*}
l\frac{\partial N}{\partial S} +\frac{\partial N}{\partial t} + \frac{\partial N}{\partial l} \left(\cot(\alpha_1)+\cot(\alpha_2)\right)+
\lambda N \left(G_1(\alpha_1)+G_2(\alpha_2)\right)-
\\
F(\alpha_1)\cdot\int\limits_{0}^{\alpha_1}
N(S,l,\varphi,\alpha_2,t)\cdot\frac{\sin(\alpha_1-\varphi)}{\sin(\varphi)}\mathrm{d} \varphi-
\\
F(\alpha_2)\cdot\int\limits_{0}^{\alpha_2}
N(S,l,\alpha_1,\varphi,t)\cdot\frac{\sin(\alpha_2-\varphi)}{\sin(\varphi)}\ \mathrm{d} \varphi.
\end{multline*}

where
$$
G_1(\alpha)=\int\limits_{\alpha}^{\pi}
F(\varphi)\cdot\frac{\sin(\varphi-\alpha)}{\sin(\alpha)}\ \mathrm{d} \varphi \
\text{and} \
   G_2(\alpha)=\int\limits_{\alpha}^{\pi}
F(\pi-\varphi)\cdot\frac{\sin(\varphi-\alpha)}{\sin(\alpha)}\ \mathrm{d} \varphi.
$$
\end{theorem}

\begin{theorem} {\label{thm2}}
The transition function
$\alpha\xrightarrow[l,P]{t}\alpha'$,
$\alpha\xrightarrow[l,S]{t}\alpha'$ and
$\alpha\xrightarrow[P,S,l]{t}\alpha'$
can be formulate as following:
\begin{equation*}
\left(l\frac{\partial}{\partial S}+\frac{\partial}{\partial t}+Q_1^F\right)\cdot\alpha\xrightarrow[l,S]{t}\alpha'=0,
\end{equation*}

\begin{equation*}
\left(\frac{1}{\sin(\alpha')}\frac{\partial}{\partial P}+\frac{\partial}{\partial t}+Q_1^F\right)\cdot\alpha\xrightarrow[l,P]{t}\alpha'=0,
\end{equation*}

\begin{equation*}
\left(l\frac{\partial}{\partial S}+\frac{1}{\sin(\alpha')}\frac{\partial}{\partial P}+\frac{\partial}{\partial t}+Q_1^F\right)\cdot\alpha\xrightarrow[l,S]{t}\alpha'=0,
\end{equation*}

where

\begin{multline*}
Q_1^F\left(\Re(*,*,l,\alpha,\alpha')\right)=\cot(\alpha)\cdot\frac{\partial\Re}{\partial l}+G_1(\alpha')\cdot\Re - \\
F(\alpha')\cdot\int\limits_{\alpha}^{\alpha'}
\frac{\sin(\alpha'-\varphi)}{\sin(\varphi)}\mathrm{d} \varphi \cdot
\Re(*,*,l,\alpha,\varphi).
\end{multline*}

And the distribution of line in Theorem \ref{thm1} can be expressed as following.
Where does the area distribution function look like:

\begin{equation*}
\mathcal{N}(S,t)=
\frac{1}{C^2}\cdot\int\limits_{\mathcal{D}_5}
\alpha_1^O\xrightarrow[l,1]{t}\alpha_1\cdot
\alpha_2^O\xrightarrow[l,2]{t}\alpha_2 \
\mathrm{d} l_1 \mathrm{d} S_1
\mathrm{d}\alpha_1^O\mathrm{d}\alpha_2^O,
\end{equation*}

where $\mathcal{D}_5$ is set by conditions:
$$l_1+l_2=l, \ S_1+S_2 =S, \ \alpha_1^O,\alpha_2^O>0, \ \alpha_1^O+\alpha_2^O <\pi .$$

Similar formulas can be written for perimeters:
\begin{equation*}
\mathcal{N}(P,t)=
\frac{1}{C^2}\cdot\int\limits_{\mathcal{D}_6}
\alpha_1^O\xrightarrow[l,1]{t}\alpha_1\cdot
\alpha_2^O\xrightarrow[l,2]{t}\alpha_2 \
\mathrm{d} l_1 \mathrm{d} P_1
\mathrm{d}\alpha_1^O\mathrm{d}\alpha_2^O
\end{equation*}

where $\mathcal{D}_6$ is set by conditions:
$$l_1+l_2=l, \ P_1+P_2 =P, \  \alpha_1^O,\alpha_2^O>0,  \ \alpha_1^O+\alpha_2^O <\pi $$
and for joint distribution by area and perimeter:

\begin{equation*}
\mathcal{N}(S,P,t)=
\frac{1}{C^2} \cdot\int\limits_{\mathcal{D}_7}
\alpha_1^O\xrightarrow[l,S]{t}\alpha_1\cdot
\alpha_2^O\xrightarrow[l,S]{t}\alpha_2 \
\mathrm{d} l_1 \mathrm{d}P_1
\mathrm{d}S_1
\mathrm{d}\alpha_1^O\mathrm{d}\alpha_2^O,
\end{equation*}
where $\mathcal{D}_7$ is set by conditions:
$$l_1+l_2=l,\ P_1+P_2=P,\ S_1+S_2=S,\ \alpha_1^O,\alpha_2^O>0, \  \alpha_1^O+\alpha_2^O<\pi,$$
\end{theorem}

\begin{theorem} {\label{thm3}}
The distribution function of area, perimeter and joint distribution of area and perimeter will be represented through of transitions.
More precisely,
\begin{equation*}
\mathcal{N}(S)=\frac{1}{Q_S}\cdot
\int\limits_{t>0}
\frac{\mathcal{N}(S,O,t)}{t} \   \mathrm{d}t,
\end{equation*}
such that
\begin{equation*}
Q_S=\iint\limits_{t>0} \frac{\mathcal{N}(S,O,t)}{t} \    \mathrm{d} t \mathrm{d} S.
\end{equation*}

Similarly

\begin{equation*}
\mathcal{N}(P)=\frac{1}{Q_P}\cdot
\int\limits_{t>0}
\frac{\mathcal{N}(P,O,t)}{t} \   \mathrm{d}t,
\end{equation*}

such that

 \begin{equation*}
Q_P=\iint\limits_{t>0}\frac{\mathcal{N}(P,O,t)}{t} \   \mathrm{d}t \mathrm{d} P,
\end{equation*}

and

\begin{equation*}
\mathcal{N}(S,P)=\frac{1}{Q_{S,P}}\cdot
\int\limits_{t>0}
\frac{\mathcal{N}(S,P,O,t)}{t} \
  \mathrm{d}t,
\end{equation*}
such that
\begin{equation*}
Q_{S,P}=\iiint\limits_{t>0}\frac{\mathcal{N}(S,P,O,t)}{t} \   \mathrm{d} t \mathrm{d} S \mathrm{d} P.
\end{equation*}

\end{theorem}

\begin{corollary}
For Voronoi diagram, these  results are satisfied.
\end{corollary}

Consider a Poisson point process on a plane. With each such point, we associate its area of attraction, i.e., the set of points for which this point is closer than any other marked one. The problem is to find distributions. In fact, the distributions are the same as for the random field of lines \cite{2}. Let us explain why. In the space of lines there is a measure that is invariant under the motions of the plane
A straight line can be given by an angle $\alpha$ and a distance $p$ to a fixed point $O$, and this measure is $\mathrm{d}\alpha \mathrm{d}p$.

If we consider poisson point process on a plane and fix point $0$ from that poisson set, and choose any other poison point $P$ then the propability density corresponding line dividing plane on areas closer to $O$ and close to $P$ would have same density proportional  $\mathrm{d}\alpha \mathrm{d}p$ where $p$ is the distant from $P$ to $O$.

When we study a section of a moving straight line $L$ of a region $P$ (a Voronoi polygon), as $O$ we take its point of the Poisson process corresponding to $P$ i.e., such that any other point in the process is further from further from any other point than $O$.

Then we see that the transition probabilities are the same as in the case of splitting by straight lines (however similar)\cite{2}.


\section{Preliminaries}   \label{ScPrel}
We make use of common notions in stochastic geometry and usual definitions
from differential equations.

\subsection{Convention and notations}    \label{ScConvNot}
The set of rational numbers, real numbers, and complex numbers is denoted by $\mathbb{Q}$, $\mathbb{R}$, and $\mathbb{C}$, respectively.
In addition, if $z \in \mathbb{C}$, then real and imaginary part of $z$ is denoted by
$\Re(z)$ and $\Im(z)$, respectively.

Suppose that $A,B$ be two points in the plane.
Then
we denote the segment that connects the two points $A$ and $B$ by $[A,B]$, and its length is denoted by
$|AB|$.

Let $f, g : \mathbb{R} \to \mathbb{R}$ be two functions.
We
implicitly assume that the following notion is defined for
$x \to +\infty$.
We write $f = \mathcal{O}(g)$ ( or $f \in \mathcal{O}(g)$), if there exist two constants $x_0$
and $D > 0$ such that, for all $x\geqslant x_0$,
$|f(x)| \leqslant D |g(x)|$.

\subsection{Poisson line process}   \label{SbScPoisonLine}
We can interpret a line process as a point process.
\begin{definition}
A {\it Poisson line process} with intensity $\lambda$ is parametrized by a Poisson point
process with intensity $\lambda$ on $[0, \pi) \times \mathbb{R}$.
\end{definition}

\begin{figure}[ht]
    \centering
    \includegraphics[width=0.7\textwidth]{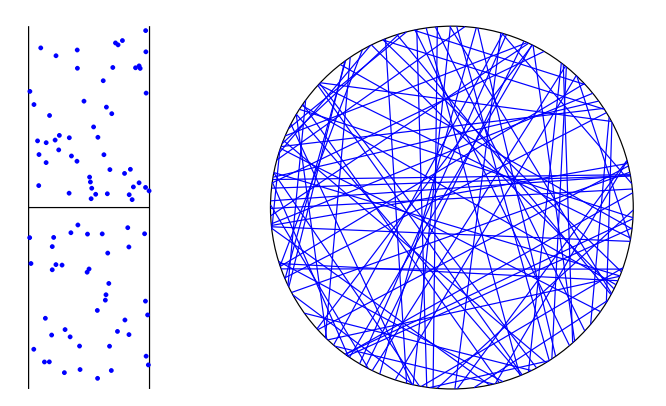}
    \caption{Poisson line process}
\end{figure}

By designating a direction to each line, undirected Poisson line processes can be transformed into Poisson directed line processes. Consequently, a natural parameter space would be $\lambda$ on $[0, 2\pi) \times \mathbb{R}$.

Poisson polygons are delimited by Poisson lines.

There are two unique approaches to express the distribution associated with a Poisson polygon attribute (area, perimeter, number of vertices):

\begin{itemize}
    \item
   Assign the same weight to every polygon in a single realisation. The obtained statistical properties are those of a hypothetical polygon known as the typical polygon.

    \item
Consider the statistical characteristics of the polygon containing a fixed point, such as the origin (fundamental polygon).
\end{itemize}

\section{Random plane partitions. Kinetic equations.}  \label{ScKinetic}

 Consider the following problems:

{\bf 1.}
On the plane, a random set of straight lines is given, all shifts are equally likely, and the law of direction distribution has the form $F(\varphi)$.

What is the area distribution of the partition's components?
Joint distributions are of interest.
We shall focus on the \textit{Poisson field} of straight lines in particular.

For instance, book \cite{2} contains a concise exposition of this topic.

{\bf 2.}
A random set of points is marked on the plane.
Each point $A$ of this set is associated with its area of attraction, i.e., the set of points of the plane for which the point $A$ is the closest of the marked ones.
The regions of attraction are convex polygons, and the same questions arise for them.

This partition is called the \textit{Voronoi diagram}, a strict statement of this problem is contained in book \cite{2}.

The following auxiliary ideas will be used in the future:

Due to the independence of the movement of the section's ends, we can consider only one end. The following arguments alleviate the difficulties connected with the consideration of improper integrals:
the number of partitioning parts larger than $d$ decreases as $e^{-\lambda d};$ the number of parts of the area greater than $S$ decreases as $e^{-\lambda S}$ with growth, respectively, $d$ and $S.$

We will begin with problem {\bf 1} for a \textit{Poisson set} of straight lines, in which not only shifts, but also all angles of inclination of straight lines are equally probable.

Let $L$ be a secant straight line moving with unit velocity along the plane and at the same time not changing its direction.

Consider the partitioning polygons that $L$ intersects and, above all, their parts lying \textit{under} direct $L$.
Note that the distribution law for such parts differs from the natural law (in relation to the unit area) and the transition to the natural law requires special consideration.

Consider the $M$ partitioning polygon.

\begin{figure}[ht]
    \centering
    \includegraphics[width=1.0\textwidth]{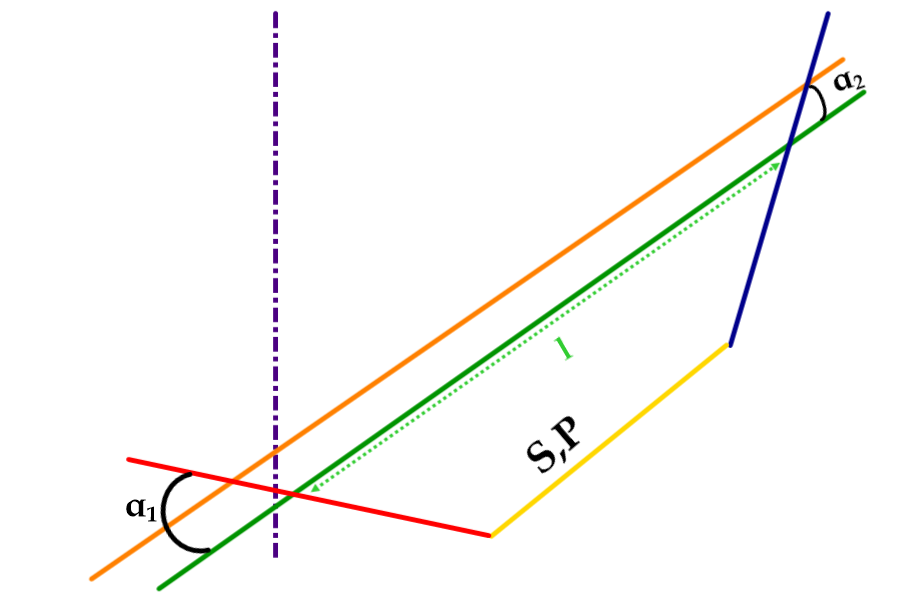}
    \caption{Changing parameters when shifting the secant line}
  \label{fig1}
\end{figure}

The traversed area is denoted by $S$, the section length of the polygon $M$ is denoted by $l$, and the angles produced by $L$ from the side of $M$ are denoted by $\alpha_1$ and $\alpha_2$
(see Figure \ref{fig1}).

Suppose that $\alpha$ is the old angle and $\alpha'$ is the new angle.
Angle change with probability
$$p\mathrm{d} t=\mathrm{d} p=\lambda \mathrm{d}t\frac{\sin(\alpha'-\alpha)}{\sin(\alpha)}.$$

In the previous figure $l$, $S$ and $P$  represent section length, area covered and passed perimeter.
Also, $\alpha_1,\alpha_2$  represent angles between cross-section and lateral straight lines, respectively.
Therefore, $\mathrm{d}l$ represents increment of the section length, then
$$\mathrm{d}l = \mathrm{d}t \left( \cot (\alpha_1) + \cot (\alpha_2) \right).$$
Therefore
$$\mathrm{d}S = l \mathrm{d}t,$$
and
$$\mathrm{d} P=\mathrm{d} t\left(\frac{1}{\sin(\alpha_1)}+\frac{1}{\sin(\alpha_2)}\right).$$

\medskip

When the straight line $L$ is shifted by the distance $\mathrm{d} t$, the following events are possible:

\begin{enumerate}
    \item[{\bf{1.}}]
    A shifted straight line $L'$ will not intersect any other side of $M$. In this case,
\begin{equation}\label{e1}
S'=S+l \mathrm{d}t;\  l'=l+\mathrm{d}t\left(\cot(\alpha_1)+\cot(\alpha_2)\right);
\ \alpha_1'=\alpha_1,\ \alpha_2'=\alpha_2.
\end{equation}

\item[{\bf{2.}}]
 The straight line $L'$ will intersect one new side of $M$.
 In this case,
 $\mathrm{d}s=l \mathrm{d}t$, and one of the corners of $\alpha_1$ or $\alpha_2$ will change.

Due to symmetry, it is sufficient to consider the case when $\alpha_1$ changes.
The new side has a greater angle of deflection $\alpha_1'$ than the angle of inclination of the old $\alpha_1$.

Let us find the probability density $Q_{\alpha}^{\alpha'}$ of the last transition. We neglect the value of $ \mathrm{d}t$ (i.e., the probability that a section with two straight lines will occur at once for $\mathrm{d}t$).

Consider the segment $[A,A']$. Its length is $1/\sin(\alpha_1).$ Probability of transition $\alpha\rightarrow\beta$, where $\alpha'<\beta<\alpha'+ \mathrm{d}\alpha'$ there is a probability that the segment $[A,A']$ will be crossed by a straight line with such an angle of inclination. This probability is equal to the projection of the segment $[A,A']$ onto the normal to the line forming with $L$ the angle $\alpha'$ ($\mathrm{d}\alpha'$  small) multiplied by $\lambda\mathrm{d}\alpha'$, where the parameter $\lambda$ characterizes the intensity of the {\it Poisson process} of straight lines.

The corresponding projection $[A,A']$ be equal to $$|AA'|\sin(\alpha'-\alpha) \ \ \text{or} \ \  \frac{\mathrm{d} t\sin(\alpha'-\alpha)}{\sin(\alpha)}.$$

Hence, the desired probability density is
$$\frac{\lambda \mathrm{d}t\sin(\alpha'-\alpha)}{\sin(\alpha)}+\mathcal{O}(\mathrm{d} t^2).$$

\item[{\bf{3.}}]
A straight line intersects two new sides of $M$.

The probability of this event is on the order of $\mathcal{O}(\mathrm{d}t^2)$, so we neglect it.
Since the number of sides of the polygon is finite, in case 2 we can assume that $\mathrm{d} l=\mathrm{d} t\left(\cot(\alpha_1)+\cot(\alpha_2)\right),$ since the error obtained for an integer polygon has the order $\mathrm{d} t\cdot N,$ where $N$ is the number of sides of the polygon.
Aiming $\mathrm{d}t$ to $0$ we get an exact result.

Through $N(S,l,\alpha_1,\alpha_2,t)$ we denote the density of the number of traversed parts adjacent to $L$ with the traversed area $S$, the length of the trace $l$, the angles $\alpha_1,\alpha_2$ at the ends of the trace, and the lowest vertex of which is located at a distance $t$ from $L$ ( see Figure \ref{fig1}).

Let us begin by defining $N$ strictly in terms of the limit:

\begin{equation}\label{e2}
N(S,l,\alpha_1,\alpha_2,t)=\lim \frac{W_{S,l,\alpha_1,\alpha_2,t}^{\Delta S,\Delta l,\Delta \alpha_1,\Delta \alpha_2,\Delta t}}{\Delta S\cdot\Delta l\cdot\Delta \alpha_1\cdot\Delta \alpha_2\cdot\Delta t},
\end{equation}
where $W$ is the relative number of parts (in unit of section length) whose parameters are enclosed within
$$[S,S+\Delta S],[l,l+\Delta l],[\alpha_1,\alpha_1+\Delta \alpha_1],[\alpha_2,\alpha_2+\Delta \alpha_2] \ \text{and} \ [t,t+\Delta t].$$
\end{enumerate}

\medskip

The parameter $t$ is introduced to account for the effect that when the secant line $L$ moves on $\mathrm{d}t$, as mentioned above, the events described in cases 1,2 and 3 occur. It is also needed for recalculation from distributions along the secant line to distributions among whole plane (see section \ref{SctoUsualDistr}). That is because each polygon during moving process will be counted with weight $t$.

Let us write the balance ratio down:

\begin{multline}\label{e3}
N(S,l,\alpha_1,\alpha_2,t)=N\left(S-l\mathrm{d} t,l- \mathrm{d}t \left(\cot(\alpha_1)+\cot(\alpha_2)\right),\alpha_1,\alpha_2,t-\mathrm{d}t\right)-\\
N(S,l,\alpha_1,\alpha_2,t) \cdot\lambda \mathrm{d} t
\left( \int\limits_{\alpha_1}^{\pi}\frac{\sin(\varphi-\alpha_1)}{\sin(\alpha_1)} \mathrm{d} \varphi+
\int\limits_{\alpha_2}^{\pi}\frac{\sin(\varphi-\alpha_2)}{\sin(\alpha_2)} \mathrm{d} \varphi  \right)+\\
\lambda \mathrm{d} t\left(\int\limits_{0}^{\alpha_1}N(S,l,\varphi,\alpha_2,t) \frac{\sin(\alpha_1-\varphi)}{\sin(\varphi)} \mathrm{d} \varphi+
\int\limits_{0}^{\alpha_2}N(S,l,\alpha_1,\varphi,t) \frac{\sin(\alpha_2-\varphi)}{\sin(\varphi)} \mathrm{d} \varphi \right).
\end{multline}

The first term of the first part corresponds to the absence of intersections for $\mathrm{d}t$, the second term corresponds to the departure from the state $(S,l,\alpha_1,\alpha_2,t)$, and the third to the arrival in this state.

\medskip

Rewrite equation \eqref{e3} as a kinetic equation:

\begin{multline}\label{e4}
l\frac{\partial N}{\partial S} + \left(\cot (\alpha_1) +\cot (\alpha_2)\right)\frac{\partial N}{\partial l} +\frac{\partial N}{\partial t} +  \lambda N \left( \tan (\frac{\alpha_1}{2}) + \tan (\frac{\alpha_2}{2})\right )-\\
\lambda  \cdot \int\limits_0^{\alpha_1} \frac{N(S,l,\varphi,\alpha_2,t) \sin (\alpha_1-\varphi)}{\sin(\varphi)} \mathrm{d} \varphi -
\lambda  \cdot \int\limits_0^{\alpha_2} \frac{N(S,l,\alpha_1,\varphi,t) \sin (\alpha_2-\varphi)}{\sin(\varphi)} \mathrm{d} \varphi =0.
\end{multline}

\medskip

Similarly, you can write out the kinetic equation of the perimeters.
Let $p$ pass the perimeter:
\begin{multline}\label{e5}
\left(\frac{1}{\cos (\alpha_1)} + \frac{1}{\cos (\alpha_2)}\right)\frac{\partial N}{\partial P} + \frac{\partial N}{\partial t}
+\frac{\partial N}{\partial l}\left(\cot (\alpha_1) +\cot (\alpha_2)\right)+ \lambda N \left( \tan(\frac{\alpha_1}{2}) + \tan(\frac{\alpha_2}{2})  \right) -\\
\lambda  \cdot \int\limits_0^{\alpha_1} \frac{N(P,l,\varphi,\alpha_2,t) \sin (\alpha_1-\varphi)}{\sin(\varphi)} \mathrm{d} \varphi -
\lambda  \cdot \int\limits_0^{\alpha_2} \frac{N(P,l,\alpha_1,\varphi,t) \sin (\alpha_2-\varphi)}{\sin(\varphi)} \mathrm{d} \varphi =0.
\end{multline}
As a remark we need to elaborate, the integrals are conversion duo to $N \to \infty$, if $\alpha \to 0$.

Now, let us write the {\it kinetic equation} for the case of an arbitrarily chosen law $F(\varphi)$ governing the distribution of inclination angles.
Assume that the angle of inclination of the line $L$ is zero.

Let $\mathrm{d}F$ be the expectation of the number of intersection points of a unit segment normal to a straight line with an angle of deviation $\varphi$ with a straight line whose angle of deviation is enclosed between $\varphi$ and $\varphi+ \mathrm{d} \varphi$.
We define the function $F(\varphi)$ as the limit:

\begin{equation}\label{e6}
F(\varphi)=\lim_{\mathrm{d} \varphi\rightarrow 0}\frac{\mathrm{d} F}{\mathrm{d}\varphi}.
\end{equation}

For the left end, the probability density of the transition $\alpha_1\rightarrow \alpha_1'$ for {\it time} $\mathrm{d}t$ be equal to $$\mathrm{d} t\cdot F(\alpha_1')\cdot \frac{\sin(\alpha_1'-\alpha_1)}{\sin(\alpha_1)},$$
and for the right end, the probability density of the transition $\alpha_2\rightarrow \alpha_2'$ for {\it time} $\mathrm{d}t$ be equal to $$\mathrm{d} t\cdot F(\alpha_2')\cdot \frac{\sin(\alpha_2'-\alpha_2)}{\sin(\alpha_2)}.$$

The kinetic equation of areas can be written as:
\begin{multline}\label{e7}
l\frac{\partial N}{\partial S} +\frac{\partial N}{\partial t} + \frac{\partial N}{\partial l} \left(\cot(\alpha_1)+\cot(\alpha_2)\right)+
\lambda N \left(G_1(\alpha_1)+G_2(\alpha_2)\right)-
\\
F(\alpha_1)\cdot\int\limits_{0}^{\alpha_1}
N(S,l,\varphi,\alpha_2,t)\cdot\frac{\sin(\alpha_1-\varphi)}{\sin(\varphi)}\mathrm{d} \varphi-
\\
F(\alpha_2)\cdot\int\limits_{0}^{\alpha_2}
N(S,l,\alpha_1,\varphi,t)\cdot\frac{\sin(\alpha_2-\varphi)}{\sin(\varphi)}\ \mathrm{d} \varphi,
\end{multline}

where

\begin{equation}\label{e8}
G_1(\alpha)=\int\limits_{\alpha}^{\pi}
F(\varphi)\cdot\frac{\sin(\varphi-\alpha)}{\sin(\alpha)}\ \mathrm{d} \varphi,
\end{equation}

and

\begin{equation*}
   G_2(\alpha)=\int\limits_{\alpha}^{\pi}
F(\pi-\varphi)\cdot\frac{\sin(\varphi-\alpha)}{\sin(\alpha)}\ \mathrm{d} \varphi.
\end{equation*}

\medskip

Consider the following equation for the perimeters:
\begin{multline}\label{e9}
\left(\frac{1}{\sin (\alpha_1) }+ \frac{1}{\sin (\alpha_2) }  \right)\frac{\partial N}{\partial P} + \frac{\partial N}{\partial t}
+\frac{\partial N}{\partial l}\left(\cot (\alpha_1) +\cot (\alpha_2)\right)+ N\cdot \left(G_1(\alpha_1)+G_2(\alpha_2)\right)
 -\\
F(\alpha_1)  \cdot\int\limits_0^{\alpha_1} N(P,l,\varphi,\alpha_2,t)\cdot\frac{\sin (\alpha_1-\varphi)}{\sin(\varphi)} \mathrm{d} \varphi -
\ F(\pi- \alpha_2) \cdot \int\limits_0^{\alpha_2} N(P,l,\alpha_1,\varphi,t) \cdot\frac{\sin (\alpha_2-\varphi)}{\sin(\varphi)} \mathrm{d} \varphi =0.
\end{multline}

Equations \eqref{e7} and \eqref{e9} can be rewritten as:

\begin{multline}\label{e7'}
l\frac{\partial N}{\partial S} + \frac{\partial N}{\partial t}
+Q_1^F[N]+Q_2^F[N]=0,\\
\frac{\partial N}{\partial P}\frac{1}{\sin(\alpha_1)} +
\frac{\partial N}{\partial P}\frac{1}{\sin(\alpha_2)} +
\frac{\partial N}{\partial t}
+Q_1^F[N]+Q_2^F[N]=0,
\end{multline}
where
\begin{multline}
Q_1^F[N]=\frac{\partial N}{\partial P}\cot (\alpha_1)+G_1(\alpha_1)\cdot N-\\
F(\alpha_1)\cdot \int\limits_0^{\alpha_1} N(*,l,\varphi,\alpha_2,t)  \frac{\sin (\alpha_1-\varphi)}{\sin(\varphi)} \mathrm{d} \varphi,\\
Q_2^F[N]=\frac{\partial N}{\partial P}\cot (\alpha_2)+G_2(\alpha_2)\cdot N-\\
F(\pi-\alpha_2)\cdot \int\limits_0^{\alpha_2} N(*,l,\varphi,\alpha_1,t) \frac{\sin (\alpha_2-\varphi)}{\sin(\varphi)} \mathrm{d} \varphi,
\end{multline}
where $*$ means the area in the first case and the perimeter in the second, it acts as a parameter in all operators $Q_i^F$.

\medskip

Similarly, a joint equation is written for the area and perimeter:

\begin{multline}\label{e12}
l\frac{\partial N}{\partial S}
+\frac{\partial N}{\partial P}\frac{1}{\sin(\alpha_1)}+\frac{\partial N}{\partial P}\frac{1}{\sin(\alpha_2)}+\frac{\partial N}{\partial t}
+Q_1^F[N]+Q_2^F[N]=0.
\end{multline}

The $Q_i^F$ operators commute.
The independence of the local movement of the ends of segments is reflected in the form of equations \eqref{e7}, \eqref{e9} and \eqref{e7'}:
the operators associated with $f$ and $g$ are independent.

Let
$$\widetilde Q_i^F=Q_i^F+\frac{\partial N}{\partial P}\frac{1}{\sin(\alpha_i)}.$$

Then equations \eqref{e9} and \eqref{e12} will be rewritten as:
\begin{equation}\label{e9''}
\frac{\partial N}{\partial t}+\widetilde Q_1^F[N]+\widetilde Q_2^F[N]=0,
\end{equation}
\begin{equation}\label{e11''}
l\frac{\partial N}{\partial S}+\frac{\partial N}{\partial t}+Q_1^F[N]+Q_2^F[N]=0,
\end{equation}
moreover, $\alpha_i$ will raise only $Q_i^F$.
This circumstance will help us in the future to get rid of one of the corners.

\subsection{Boundary conditions}
Let $t\to 0$. Then 

$$l=t\cdot\left(\cot(\alpha_1)+\cot(\alpha_2)\right), s=l\cdot t/2=\frac{t^2}{2}\left(\cot(\alpha_1)+\cot(\alpha_2)\right),$$
$$p=l+t\left(\frac{1}{\cos(\alpha_1)}+\frac{1}{\cos(\alpha_2)}\right)=
t\cdot\left(\cot(\alpha_1)+\cot(\alpha_2)+\frac{1}{\cos(\alpha_1)}+\frac{1}{\cos(\alpha_2)}\right)=$$

$$=t\cdot\left(\frac{1}{\cos(\alpha_1/2+\pi/4)}+\frac{1}{\cos(\alpha_2/2+\pi/4)}\right)$$

Angle density $\mathcal{N}(\alpha_1,\alpha_2)$ for $t=0$ and intensity of line process $\lambda$ is as follows:

$$
\mathcal{N}(\alpha_1,\alpha_2)=\lambda^2\sin(\alpha_1)\sin(\alpha_2).
$$

Thus distribution function $\mathcal{N}(s,p,t,l,\alpha_1,\alpha_2)$ asymptotically for $t\to 0$ as a generalised function tends to

\begin{multline}\label{EqBoundCond}
  \mathcal{N}(s,p,t,l,\alpha_1,\alpha_2) = \mathcal{N}(\alpha_1,\alpha_2)\cdot\delta\left(l-t\cdot\left(\cot(\alpha_1)+\cot(\alpha_2)\right)\right)\cdot \\
  \delta\left(s-\frac{t^2}{2}\left(\cot(\alpha_1)+\cot(\alpha_2)\right)\right) \cdot \delta\left(p-t\cdot\left(\frac{1}{\cos(\alpha_1/2+\pi/4)}+\frac{1}{\cos(\alpha_2/2+\pi/4)}\right)\right),
\end{multline}



where $\delta$ is a {\em Dirac $\delta$-function}.

One-sided version of the equation (\ref{EqBoundCond}) is

\begin{equation}\label{EqBoundCondOneside}
  \mathcal{N}(\alpha_1)\cdot\delta\left(l-t\cdot\left(\cot(\alpha_1)\right)\right)
\delta\left(s-\frac{t^2}{2}\left(\cot(\alpha_1)\right)\right)
\delta\left(p-t\cdot\left(\frac{1}{\cos(\alpha_1/2+\pi/4)}\right)\right)
\end{equation}

$$
\mathcal{N}(\alpha_1)=\lambda\sin(\alpha_1).
$$

\newpage
Similar equation holds for $\alpha_2$.
\begin{figure}[!ht]
    \centering
    \includegraphics[width=0.6\textwidth]{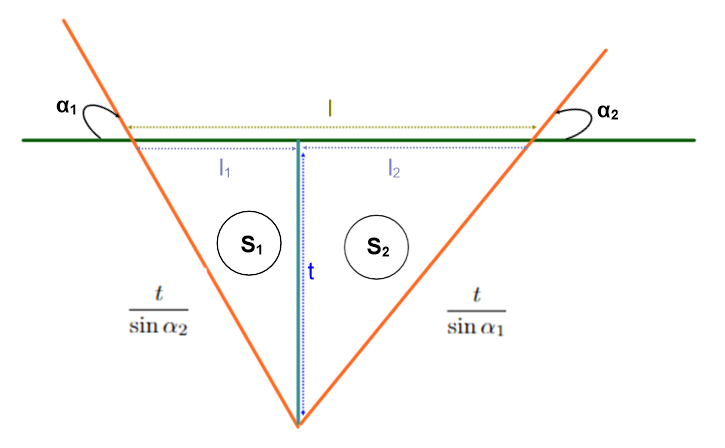}
    \caption{Boundary condition $t \to 0$}
  \label{fig2}
\end{figure}
where
$$l=l_1+l_2 =t \cdot \tan(\alpha_1) + t \cdot tan(\alpha_2),$$
$$P= P_1+P_2 = l_1 + \frac{t}{\sin(\alpha_1)} + \frac{t}{\sin(\alpha_2)}$$
and
$$S=S_1+S_2 = \frac{l_1t}{2}+ \frac{l_2t}{2} = \frac{lt}{2}.$$

\section{Local interval ends independence. Transition functions and equations simplifications} \label{ScEndsIndep}

Let us examine the behaviour of the {\it left} end with angle $\alpha_1$ (the reduction of the study of the equation \eqref{e1} and the development of the distribution function will be discussed in the section \ref{razd2-1}). We will begin by defining and designating the necessary terms.

\begin{figure}[!ht]
    \centering
    \includegraphics[width=1.0\textwidth]{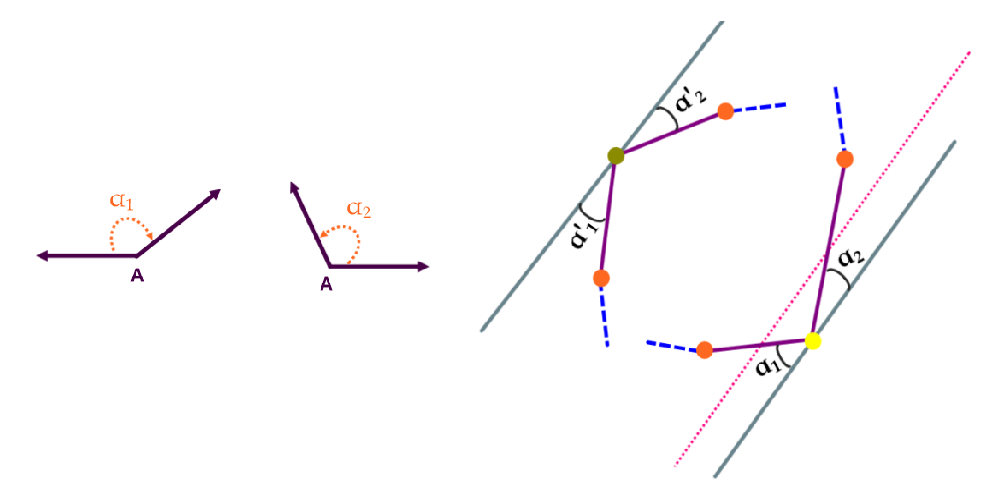}
    \caption{Convolution, independence of the movement of the ends}
  \label{fig2}
\end{figure}

\medskip

Transitions
$$\alpha_1\rightarrow\alpha_1'\  \text{and} \ \alpha_2\rightarrow\alpha_2'.$$
In this case, you can determine the increment of the perimeter $\Delta P$ and the displacement $\Delta l$.

\medskip

Consider the {\it left} the end of the section ( see Figure \ref{fig2}), its {\it state} is characterized by an angle
$\alpha_1$  when $0<\alpha_1<\pi.$
After a period of time $t$, its new state is characterized by the angle $\alpha_1'\geqslant \alpha_1.$ When switching to a new state, the end of the segment can pass the path $\Delta l$ and {\it cause} an increment of the perimeter $\Delta p.$

Let us first establish the area increment.
To do this, fix the origin of coordinates on the line $L$.

Let the point $O$ be the left end of the section.
From $O$ draw a straight line $M\perp~L.$
We will find the areas as in the Figure \ref{fig1},

$$\Delta S=\Delta S_1-\Delta S_2,$$
$\Delta S$ the area covered.

Using $\alpha\xrightarrow[\Delta S,\Delta l]{t}\alpha'$, we will denote the probability density of the transition from $\alpha$  to $\alpha'$, where the area increment is $\Delta S$ and the relocation is $\Delta l$.
Similarly, you can define the value $\alpha\xrightarrow[\Delta P,\Delta l]{t}\alpha'$  for perimeters, and for areas and perimeters simultaneously $\alpha\xrightarrow[\Delta S, \Delta P,\Delta l]{t}\alpha'.$

\medskip

Let there be two consecutive transitions $\alpha\rightarrow\alpha'$ and $\alpha'\rightarrow\alpha''$.
The first is characterized by increments of $t_1,\Delta P_1,\Delta S_1$ and $\Delta l_1,$ the second $t_2,\Delta P_2,\Delta S_2$ and $\Delta l_2$.
Then the through transition is characterized by the values
$$t_1+t_2,\Delta P_1+\Delta P_2,\Delta l_1+\Delta l_2 \ \text{and} \ \Delta S_1+t_2\Delta l_1+\Delta S_2.$$
The second term is caused by shifting the origin by $\Delta l_1.$

For the values $$\alpha\xrightarrow[\Delta S,\Delta l]{t}\alpha', \alpha\xrightarrow[\Delta P,\Delta l]{t}\alpha' \ \text{and} \ \alpha\xrightarrow[\Delta S,\Delta P,\Delta l]{t}\alpha',$$
we can write the generalized {\it Markov equation}:

\begin{equation}\label{e12a}
\alpha\xrightarrow[l,S]{t_1+t_2}\alpha'=
\int\limits_{\mathcal{D}_1}
\alpha\xrightarrow[l_1,S_1]{t_1}\beta
\cdot
\beta\xrightarrow[l_2,S_2]{t_1}\alpha' \
\mathrm{d} \beta
\mathrm{d}l_1\mathrm{d} S_1,
\end{equation}
where $\mathcal{D}_1$ is set by conditions:
$$\alpha \leqslant \beta \leqslant \alpha', \ l_1+l_2=l \ \text{and} \ S_1+S_2+l_1\cdot t_2=S,$$
\begin{equation}\label{e12'}
\alpha\xrightarrow[l,P]{t_1+t_2}\alpha'=
\int\limits_{\mathcal{D}_2}
\alpha\xrightarrow[l_1,P_1]{t_1}\beta
\cdot
\beta\xrightarrow[l_2,P_2]{t_1}\alpha' \
\mathrm{d} \beta
\mathrm{d} l_1\mathrm{d} P_1,
\end{equation}
where $\mathcal{D}_2$ is set by conditions:
$$\alpha \leqslant \beta \leqslant \alpha', \ l_1+l_2=l \ \text{and} \ P_1+P_2=P,$$
\begin{equation}\label{e12''}
\alpha\xrightarrow[P,S]{t_1+t_2}\alpha'=
\int\limits_{\mathcal{D}_3}
\alpha\xrightarrow[P_1,S_1,l_1]{t_1}\beta
\cdot
\beta\xrightarrow[P_2,S_2,l_2]{t_1}\alpha' \
\mathrm{d} \beta
\mathrm{d} l_1\mathrm{d} P_1\mathrm{d} S_1,
\end{equation}
where $\mathcal{D}_3$ is set by conditions:
$$\alpha \leqslant \beta \leqslant \alpha', \ l_1+l_2=l, \ P_1+P_2=P \ \text{and} \ S_1+S_2+l_1\cdot t_2=S.$$

Now, let us write the values' equations
$\alpha\xrightarrow[l,P]{t}\alpha'$,
$\alpha\xrightarrow[l,S]{t}\alpha'$ and
$\alpha\xrightarrow[P,S,l]{t}\alpha'$:

\begin{equation}\label{e13}
\left(l\frac{\partial}{\partial S}+\frac{\partial}{\partial t}+Q_1^F\right)\cdot\alpha\xrightarrow[l,S]{t}\alpha'=0,
\end{equation}

\begin{equation}\label{e13'}
\left(\frac{1}{\sin(\alpha')}\frac{\partial}{\partial P}+\frac{\partial}{\partial t}+Q_1^F\right)\cdot\alpha\xrightarrow[l,P]{t}\alpha'=0,
\end{equation}

\begin{equation}\label{e13''}
\left(l\frac{\partial}{\partial S}+\frac{1}{\sin(\alpha')}\frac{\partial}{\partial P}+\frac{\partial}{\partial t}+Q_1^F\right)\cdot\alpha\xrightarrow[l,S]{t}\alpha'=0,
\end{equation}
where
\begin{multline}
Q_1^F\left(\Re(*,*,l,\alpha,\alpha')\right)=\cot(\alpha)\cdot\frac{\partial\Re}{\partial l}+G_1(\alpha')\cdot\Re - \\
F(\alpha')\cdot\int\limits_{\alpha}^{\alpha'}
\frac{\sin(\alpha'-\varphi)}{\sin(\varphi)}\mathrm{d} \varphi \cdot
\Re(*,*,l,\alpha,\varphi).
\end{multline}

Equations \eqref{e13}-\eqref{e13''} are derived similarly to equations \eqref{e7},\eqref{e9} and \eqref{e7'}.

\medskip

Formal integration of equations \eqref{e13}-\eqref{e13''}:
by making the {\it Laplace transformation} in equation \eqref{e13} for $S$  and $t$, in equation \eqref{e13'} for $P, S,$ and $l$.  Moreover, in equation \eqref{e13''} for $P$ and $S$, we have:

\begin{equation}\label{e14}
\left(-\widetilde l\widetilde S-\widetilde t+ \widetilde Q_1^F\right) \cdot\widehat{\alpha\xrightarrow[\widetilde S, \widetilde l]{t}\alpha'}=H_1,
\end{equation}

\begin{equation}\label{e14'}
\left(-\frac{1}{\sin(\alpha')}\widetilde P-\widetilde t+ \widetilde Q_1^F\right) \cdot\widehat{\alpha\xrightarrow[\widetilde P, \widetilde l]{t}\alpha'}=H_2,
\end{equation}

\begin{equation}\label{e14''}
\left(-\widetilde l\widetilde S-\frac{1}{\sin(\alpha')}\widetilde P-\widetilde t+ \widetilde Q_1^F\right) \cdot \widehat{\alpha\xrightarrow[\widetilde P, \widetilde S, \widetilde l]{t}\alpha'}=H_3,
\end{equation}

where
\begin{multline*}
\widetilde Q_1^F\left(\widetilde\Re(*,*,\widetilde l,\alpha,\alpha')\right)=-\cot(\alpha')\cdot \widetilde l+
G_1(\alpha')\cdot\widetilde \Re-\\
F(\alpha')\cdot\int\limits_{\alpha}^{\alpha'}
\frac{\sin(\alpha'-\varphi)}{\sin(\varphi)} \mathrm{d}\varphi \cdot
\Re(*,*,\widetilde l,\alpha,\varphi),
\end{multline*}
and $H_1,H_2$ and $H_3$ are defined by boundary conditions.

Note that density exponentially decrease respect to perimeter and square root of arrear so Laplace transformation is correct.

The equation where the area appears cannot be integrated, yet equation \eqref{e14''} reduces to the Riccati equation.
After $\mathcal{N}(\alpha,p,t,l,\alpha')$, we denote $\alpha\xrightarrow[l,S]{t}\alpha'$  (for convenience of writing).
Then equation \eqref{e14} is rewritten as

\begin{multline}\label{e15}
\left(-\frac{1}{\sin(\alpha)}\cdot\widetilde P-\widetilde t-\cot(\alpha')\cdot\widetilde l+ G_1(\alpha')\right) \cdot N- \\
F(\alpha')\cdot
\int\limits_{\alpha}^{\alpha'}\mathcal{N}(\alpha,\widetilde p,\widetilde l,\widetilde t,\beta) \frac{\sin(\alpha'-\beta)}{\sin(\beta)} \mathrm{d} \beta
=H.
\end{multline}

Set
\begin{equation*}
\mathcal{T}(\alpha,\widetilde p,\widetilde l,\widetilde t,\alpha')=\frac{\left(-\widetilde P - \widetilde t\sin(\alpha')-\widetilde l\cos(\alpha')+G_1(\alpha')\sin(\alpha')\right)}{F(\alpha')},
\end{equation*}
and
$\widehat{H}={H}/{F(\alpha')}.$
Then equation \eqref{e15} is rewritten as:
\begin{equation}\label{e16}
\mathcal{T\widehat{N}}+\int\limits
\mathcal{\widehat{N}}(\alpha,\beta)
\sin(\alpha'-\beta)\mathrm{d} \beta=\widehat{H}.
\end{equation}

Differentiating $\alpha'$ twice and adding it together, we get:
\begin{equation*}
\mathcal{T\widehat{N}}''+\mathcal{T\widehat{N}}+
\mathcal{\widehat{N}}=\widehat{H}''+\widehat{H},
\end{equation*}
where $'$ means the derivative in respect to $\alpha'$.
From where, putting $\mathcal{Y}= \mathcal{T\widehat{N}}$ and $ \nu=\widehat{H}''+\widehat{H},$ we have:

\begin{equation*}
\mathcal{Y}''+\frac{\mathcal{Y}\cdot(\mathcal{T}+1)}{\mathcal{T}}=\nu.
\end{equation*}

Let $\mathcal{G}=(\mathcal{T}+1)/\mathcal{T}.$ Then we have the equation:
\begin{equation*}
\mathcal{Y}''+\mathcal{G}\cdot\mathcal{Y}=\nu,
\end{equation*}
which, as we know, reduces to the {\it Ricatti equation}.


\subsection{The expression of the distribution function through of transitions}
\label{razd2-1}

\medskip

\

In this part we talk about distribution of function based on transitions
$$\alpha\xrightarrow[l,P]{t}~\alpha', \alpha\xrightarrow[l,S]{t}\alpha' \ \text{and} \ \alpha\xrightarrow[S,P,l]{t}\alpha'.$$

\section{Expression of plane distriputions via distribution along the section line} \label{SctoUsualDistr}

Let us express the functions
$$\mathcal{N}(S,l,\alpha_1,\alpha_2,t), \mathcal{N}(P,l,\alpha_1,\alpha_2,t) \  \text{and} \
\mathcal{N}(S,P,l,\alpha_1,\alpha_2,t).$$
Consider an arbitrary point $O$ where the lines intersect.
Draw a straight line parallel to the secant $L$ through the point $O$.
Then the lines $l_1$ and $l_2$ form the angles $\beta_1$ and $\beta_2$ with the line $L$ (see Fig \ref{fig2}).

\medskip

The probability $\mathrm{d}p$ that $\alpha_1^O<\beta_1\alpha_1^O+\mathrm{d} \alpha_1^O$ and $\alpha_2^O<\beta_2\alpha_2^O+\mathrm{d} \alpha_2^O$ be equal to:

\begin{equation}\label{e18}
\mathrm{d} p=\frac{F(\alpha_1^O)\sin(\alpha_1^O)F(\pi-\alpha_2^O)\sin(\alpha_2^O) \mathrm{d} \alpha_1^O \mathrm{d} \alpha_2^O}{C^2},
\end{equation}
where
\begin{equation}\label{e18'}
C=\int\limits_0^{\pi} F(\alpha) \cdot \sin\alpha \ \mathrm{d} \alpha.
\end{equation}

The conditional distribution function $\mathcal{N}(S,l,\alpha_1,\alpha_2,t,\alpha_1^O,\alpha_2^O)$ has the form:
\begin{equation}\label{e19}
\mathcal{N}(S,l,\alpha_1,\alpha_2,t,\alpha_1^O,\alpha_2^O)=
\int\limits_{\mathcal{D}_4}
\alpha_1^O\xrightarrow[l,S]{t}\alpha_1\cdot
\alpha_2^O\xrightarrow[l,S]{t}\alpha_2 \
\mathrm{d} l_1
\mathrm{d}S_1,
\end{equation}
where $\mathcal{D}_4$ is set by conditions:
$$l_1+l_2=l \ \text{and} \ S_1+S_2=S.$$

Where does the area distribution function look like:

\begin{equation}\label{e19'}
\mathcal{N}(S)=
\frac{1}{C^2}\cdot\int\limits_{\mathcal{D}_5}
\alpha_1^O\xrightarrow[l,1]{t}\alpha_1\cdot
\alpha_2^O\xrightarrow[l,2]{t}\alpha_2 \
\mathrm{d} l_1 \mathrm{d} S_1
\mathrm{d}\alpha_1^O\mathrm{d}\alpha_2^O
\mathrm{d} t,
\end{equation}
where $\mathcal{D}_5$ is set by conditions:
$$l_1+l_2=l, \ S_1+S_2 =S, \ \alpha_1^O,\alpha_2^O>0, \ \alpha_1^O+\alpha_2^O <\pi \ \text{and} \ t>0.$$

Similar formulas can be written for perimeters:
\begin{equation}\label{e20}
\mathcal{N}(P)=
\frac{1}{C^2}\cdot\int\limits_{\mathcal{D}_6}
\alpha_1^O\xrightarrow[l,1]{t}\alpha_1\cdot
\alpha_2^O\xrightarrow[l,2]{t}\alpha_2 \
\mathrm{d} l_1 \mathrm{d} P_1
\mathrm{d}\alpha_1^O\mathrm{d}\alpha_2^O \mathrm{d} t,
\end{equation}
where $\mathcal{D}_6$ is set by conditions:
$$l_1+l_2=l, \ P_1+P_2 =P, \  \alpha_1^O,\alpha_2^O>0,  \ \alpha_1^O+\alpha_2^O <\pi \ \text{and} \ t>0,$$
and for joint distribution by area and perimeter:

\begin{equation}\label{e21}
\mathcal{N}(S,P)=
\frac{1}{C^2} \cdot\int\limits_{\mathcal{D}_7}
\alpha_1^O\xrightarrow[l,S]{t}\alpha_1\cdot
\alpha_2^O\xrightarrow[l,S]{t}\alpha_2 \
\mathrm{d} l_1 \mathrm{d}P_1
\mathrm{d}S_1
\mathrm{d}\alpha_1^O\mathrm{d}\alpha_2^O \mathrm{d} t,
\end{equation}
where $\mathcal{D}_7$ is set by conditions:
$$l_1+l_2=l,\ P_1+P_2=P,\ S_1+S_2=S,\ \alpha_1^O,\alpha_2^O>0, \  \alpha_1^O+\alpha_2^O<\pi \ \text{and} \ t>0.$$

\medskip

Now let us go to getting functions
$$\mathcal{N}(S,l,\alpha_1,\alpha_2,t),
\mathcal{N}(S,P,l,\alpha_1,\alpha_2,t)
\ \text{and} \
\mathcal{N}(P,l,\alpha_1,\alpha_2,t).$$

Consider the secant line $L$.
The intersection point $L$ with straight lines is characterized by the angle $\alpha$.
The angle distribution function $Q(\alpha)$  is:
$$Q(\alpha)= \frac{F(\alpha) \cdot \sin(\alpha)}{C}.$$

The intersection points themselves form a {\it Poisson set} of points on a secant line, and the angle distribution function of the various points is independent.
Find the intensity of the process $\lambda$.
It is equal to the expectation of the number of intersection points of a unit segment by $L$ or
$$\lambda=\int\limits_0^{\pi}F(\alpha) \cdot \sin\alpha \ \mathrm{d} \alpha=C.$$

Wherever the length distribution function of segments will be equal to
$$\mathcal{N}(l)= \frac{e^{-l/C}}{C}.$$
Likewise the function $\mathcal{N}(l,\alpha_1,\alpha_2)$ of distribution of traces by lengths and angles will be
$$\mathcal{N}(l,\alpha_1,\alpha_2)= \frac{e^{-1/C}\cdot F(\alpha_1)F(\pi-\alpha_2) \cdot \sin(\alpha_1)\sin(\alpha_2)}{C^3}.$$

Considering the equation backwards, we have:
\begin{equation}\label{e22}
\mathcal{N}(S,l,\alpha_1,\alpha_2,t)=
\mathcal{N}(l,\alpha_1,\alpha_2) \cdot
\int\limits_{\mathcal{D}_8}
\alpha_1^O\xrightarrow[l_1,S_1]{t}\alpha_1\cdot
\alpha_2^O\xrightarrow[l_2,S_2]{t}\alpha_2 \
\mathrm{d} l_1
\mathrm{d} S_1,
\end{equation}
where $\mathcal{D}_8$ is set by conditions:
$$l_1+l_2=l \ \text{and} \ S_1+S_2=S.$$

\begin{equation}\label{e22'}
\mathcal{N}(P,l,\alpha_1,\alpha_2,t)=
\mathcal{N}(l,\alpha_1,\alpha_2) \cdot
\int\limits_{\mathcal{D}_9}
\alpha_1^O\xrightarrow[l_1,P_1]{t}\alpha_1\cdot
\alpha_2^O\xrightarrow[l_2,P_2]{t}\alpha_2 \
\mathrm{d} l_1
\mathrm{d}P_1
\end{equation}
where $\mathcal{D}_9$ is set by conditions:
$$l_1+l_2=l \ \text{and} \ P_1+P_2=P,$$
and
\begin{multline}\label{e22''}
\mathcal{N}(S,P,l,\alpha_1,\alpha_2,t)=
\mathcal{N}(l,\alpha_1,\alpha_2)\cdot\\
\int\limits_{\mathcal{D}_{10}}
\alpha_1^O\xrightarrow[l_1,P_1,S_1]{t}\alpha_1\cdot
\alpha_2^O\xrightarrow[l_2,P_2,S_2]{t}\alpha_2 \
\mathrm{d} l_1
 \mathrm{d} P_1 \mathrm{d}S_1,
\end{multline}
where $\mathcal{D}_{10}$ is set by conditions:
$$l_1+l_2=l, \  \ P_1+P_2=P \ \text{and} \ S_1+S_2=S.$$

Let us now express $\mathcal{N}(S), \mathcal{N}(P)$ and
$\mathcal{N}(S,P)$ in terms of the corresponding section functions.
Let us use the fact that each polygon is considered with a multiplicity proportional to the length of the normal $L$.
Therefore,

\begin{equation}\label{e23}
\mathcal{N}(S)=\frac{1}{Q_S}\cdot
\int\limits_{\mathcal{D}_{11}}
\frac{\mathcal{N}(S,O,\alpha_1,\alpha_2,t)}{t} \   \mathrm{d} \alpha_1 \mathrm{d} \alpha_2 \mathrm{d}t,
\end{equation}
where $\mathcal{D}_{11}$ is set by conditions:
$$\alpha_1,\alpha_2>0, \ \alpha_1+\alpha_2<\pi \ \text{and} \ t>0,$$
such that
\begin{equation}\label{e24}
Q_S=\int\limits_{\mathcal{D}_{12}} \frac{\mathcal{N}(S,O,\alpha_1,\alpha_2,t)}{t} \   \mathrm{d} \alpha_1 \mathrm{d} \alpha_2 \mathrm{d} t\mathrm{d} S,
\end{equation}
where $\mathcal{D}_{12}$ is set by conditions:
$$S,\alpha_1,\alpha_2>0, \ \alpha_1+\alpha_2<\pi \ \text{and} \ t>0.$$

Similarly
\begin{equation}\label{e25}
\mathcal{N}(P)=\frac{1}{Q_P}\cdot
\int\limits_{\mathcal{D}_{13}}
\frac{\mathcal{N}(P,O,\alpha_1,\alpha_2,t)}{t} \  \mathrm{d} \alpha_1 \mathrm{d} \alpha_2 \mathrm{d}t,
\end{equation}
where $\mathcal{D}_{13}$ is set by conditions:
$$\alpha_1,\alpha_2>0, \ \alpha_1+\alpha_2<\pi \ \text{and} \ t>0,$$
such that
 \begin{equation}\label{e25'}
Q_P=\int\limits_{\mathcal{D}_{14}}\frac{\mathcal{N}(P,O,\alpha_1,\alpha_2,t)}{t} \  \mathrm{d} \alpha_1 \mathrm{d} \alpha_2 \mathrm{d} t\mathrm{d} P,
\end{equation}
where $\mathcal{D}_{14}$ is set by conditions:
$$P,\alpha_1,\alpha_2>0, \ \alpha_1+\alpha_2<\pi \ \text{and} \ t>0,$$
and
\begin{equation}\label{e26}
\mathcal{N}(S,P)=\frac{1}{Q_{S,P}}\cdot
\int\limits_{\mathcal{D}_{15}}
\frac{\mathcal{N}(S,P,O,\alpha_1,\alpha_2,t)}{t} \
 \mathrm{d} \alpha_1 \mathrm{d}\alpha_2 \mathrm{d}t,
\end{equation}
where $\mathcal{D}_{15}$ is set by conditions:
$$\alpha_1,\alpha_2>0, \ \alpha_1+\alpha_2<\pi \ \text{and} \ t>0,$$
such that
\begin{equation}\label{e26'}
Q_{S,P}=\int\limits_{\mathcal{D}_{16}}\frac{\mathcal{N}(S,P,O,\alpha_1,\alpha_2,t)}{t} \  \mathrm{d} \alpha_1 \mathrm{d} \alpha_2 \mathrm{d} t\mathrm{d} S\mathrm{d} P,
\end{equation}
where $\mathcal{D}_{16}$ is set by conditions:
$$S,P,\alpha_1,\alpha_2>0, \ \alpha_1+\alpha_2<\pi \ \text{and} \ t>0.$$


\section{Justification questions}
The preceding paragraph discussed the derivation of the primary equation and the shift from secant line distributions to the natural distribution law on a physical level of rigour.

This section is devoted to the issues of justification.
We will show how to bring the reasoning to a mathematical level of rigor by approximating the distribution function for a Poisson set of straight lines by the case of $n$ systems, where $n \to \infty$.

Finally, although a strict definition of the concept of a Poisson set is contained in the literature (see, for example, \cite{2,3}), however, we will choose a methodically different way to define Poisson sets of lines using the zero–one law.


\subsection{Zero–one law}
We will show how this law can be used to justify concepts related to a random environment.

To illustrate the idea, let us first consider a well-known model example.
Let $\varepsilon_i$  for all $i=1,\dots,\infty$ be a sequence of independent random variables that take the values $0$ or $1$.
Let $a_i$ be a sequence of real numbers.
What can we say about the $\sum a_i\varepsilon_i$ series?
Without knowing anything about the $a_i$ sequence, we can immediately say that it either converges or diverges with a probability of $1$.
(In fact, it converges almost surely if $\sum a_i^2<\infty$, otherwise it diverges whenever $\sum a_i^2=\infty$).

The reason for convergence with probability $0$ or $1$ is that the convergence property does not depend on any finite set of values $\varepsilon_i$.
Thus the set $M$ of convergence is independent of any set of values $\varepsilon_i$, and hence of the {\it Borel field} generated by them.
But the convergence set is defined in terms of $\varepsilon_i$, and thus belongs to the {\it Borel field} generated by these quantities.
Therefore, the set $M$ is independent of itself, so $M$ has a measure of $0$ or $1$.

A stricter justification of the zero–one law, is contained, for example, in the book \cite{6}.

If a random environment is given by an infinite set of random variables, and its asymptotic properties do not depend on any finite set, then the asymptotic statements with probability $0$ or $1$ are executed.

This remark is often used tacitly: a rock sample is referred to as a random environment, and its statistical features are thought to dictate how the sample behaves.
However, if the zero–one law was incorrect, then the environment's effective characteristics could not always be stated: there would be a fundamental difference in the properties of two identical pieces of rock; if this characteristic was set with a probability greater than $0$ or $1$, there would be a fundamental difference in the definition of effective characteristics; however, when calculating a particular characteristic, its existence is always assumed.

Let us now proceed to consider the Poisson set of lines.
The straight line is defined by the angle $\varphi$ and the distance $r$ to the origin.
In this case, $0<\varphi<\pi$ and $r > 0$.

A Poisson point process in the parameter space $(r,\varphi)$ defines a Poisson linear process on the plane.
A Poisson point process can be defined using a countable set of independent random variables.
Partitioning the parameter space into rectangles:
$K<r<K+1$, and consider $n_K$ random variables that express the number of points in the $K^{\text{th}}$ rectangle, and $T_{i,l}$ for $ i,l=1,\dots, \infty$ coordinates of the $(l-1)^{\text{th}}$  random point in the $i^{\text{th}}$  rectangle.
Then
$n_K$ takes with probability $({\alpha^n}/n!) \cdot \exp(-n \alpha)$  the value $n$. Also,  $T_{i,l}^r$  and $T_{i,l}^\varphi$  values evenly distributed on the $i^{\text{th}}$ rectangle: if $l > n$, the point they set is not taken into account.

It is clear that changing the value of a finite set of variables only affects the finite set of lines of a random partition, and therefore does not affect the asymptotic properties of the partition, so that each asymptotic statement is satisfied with a probability of $0$ or $1$.


\subsection{Approximation of the Poisson process by the case of $n$ systems.}
We will approximate the Poisson set of direct $n$ equally spaced direct systems, the direct system will be in general position, i.e., if an arbitrary pair of systems and using affine transformations of the plane take it to the lattice of unit squares in the coordinate system associated with this grating,
the tangents of the angles of the other systems are irrational numbers linearly independent over $\mathbb{Q}$, it is just the formula for the number of partition parts satisfying this property from the concepts of
ergodic distribution.
In this case, the probability of the number of points where the segment intersects with straight lines will tend to the corresponding value for the Poisson law.
The shape of the $n$-gon can be set by $n$ parameters and the corresponding distribution functions will tend to what is needed.

On the other hand, when conducting a reasoning with $N \mathrm{d}t^{-2},$ all the reasoning related to the derivation of the difference equation becomes strict.
The aspiration of $N$ to infinity completes the justification.

In other points, the justification is constructed similarly.

\section{Concluding remarks}

Richard Miles discovered partial results relevant to the topic we posed in this work in 1972.
He discussed the mean values of random division \cite{7}, but we formulated density functions and solved it entirely.

In conclusion, there are further intriguing challenges of this type, and we plan to investigate them in the future. For instance, extending the same problem to various structures, such as hyperbolic geometry, spherical geometry, etc.


\begin{thebibliography}{99}

\newcommand{\by}[1]{{\it{#1}~}}
\newcommand{\paper}[1]{{\rm{#1}. }}

\def \journ{}
\def \jour{}
\def \book{}
\def \yr{}
\def \vol{Vol}
\def \no{№}
\def \pages{p}

\bibitem{13}
\by{A. Kanel-Belov.}
\paper{Geometric properties of block media},
Dep. in VINITI. No.272-B91, VINITI, Moscow, 1991. ( in Russian)


\bibitem{14}
\by{A. Kanel-Belov.}
\paper{On random partitions},
Dep. in VINITI. No.273-B91, VINITI, Moscow, 1991. ( in Russian)

\bibitem{11}
\by{A. Kanel-Belov.}
\paper{Statistical geometry and equilibrium of block arrays},
PhD dissertation in Mathematics and Physics,
Moscow State University, Moscow, 1991. ( in Russian)

\bibitem{15}
\by{A. Kanel-Belov, M. Golafshan, S. Malev, R.Yavich.
}
\paper{ABOUT RANDOM SPLITTING OF THE PLANE},
Crimean Autumn Mathematical School-Symposium (KROMSH),
p:
294-295
2020.



\bibitem{7}
\by{R. Miles.}
\paper{The random division of space},
Advances in Applied Probability, 1972.


\bibitem{miles}
\by{R. Miles.}
\paper{Poisson flats in euclidean spaces},
 Adv. Appl. Prob., Vol. 1, pp. 211-237.


\bibitem{1}
  \by{P. Moran, M. Kendall.}
  \paper{Geometrical Probability},
 Andesite Press, 2015.

\bibitem{2}
  \by{L. Santalo.}
  \paper{Integral Geometry and Geometric Probability},
  Cambridge Mathematical Library, 1976.

\bibitem{3}
  \by{G. Matheron.}
  \paper{Random sets and integral geometry},
 Wiley series in probability and mathematical statistics, 1974.

  \bibitem{5}
  \by{R. Ambartsumyan, Y. Mekke, D.Shtoyan.}
  \paper{Introduction to Stochastic Geometry},
 Nauka, 1989. ( in Russian)

 \bibitem{4}
  \by{N. Anoshchenko.}
  \paper{Geometric analysis of fracturing and blockiness of facing stone deposits},
   Mathematical methods and applications. The Fourth Mathematical Symposium.
Moscow State University, 1983. ( in Russian)




 \bibitem{6}
  \by{J. Lamperti.}
  \paper{Probability. A survey of the mathematical theory}
Mathematics Monograph Series. New York-Amsterdam: W.A. Benjamin, Inc. X, 150 p. (1966).


\end{thebibliography}
\end{document}